\documentclass[12pt,leqno]{article}
\usepackage{amsmath, amsthm, amsfonts, amssymb, color}
\usepackage{mathrsfs}
\newtheorem{thm}{Theorem}[section]

\newtheorem{prp}[thm]{Proposition}

\theoremstyle{definition}

\newcommand{\scr}[1]{\mathscr #1}
\definecolor{wco}{rgb}{0.5,0.2,0.3}

\numberwithin{equation}{section} \theoremstyle{remark}

\newcommand{\ua}{\uparrow}

\title{
{\bf Asymptotics of  Sample  Entropy Production Rate for  Stochastic Differential Equations\footnote{F.Y. Wang is supported in part by NNSFC(11131003, 11431014), the 985 project and the Laboratory of Mathematical and  Complex Systems.
J. Xiong was  supported by Macao Science and Technology Fund FDCT 076/2012/A3 and Multi-Year Research Grants of the University of Macau Nos. MYRG2014-00015-FST and MYRG2014-00034-FST. L. Xu is supported by the grant Science and Technology Development Fund, Macao S.A.R FDCT  049/2014/A1 and the grant MYRG2015-00021-FST.
}}}

\author{
{\bf Feng-Yu Wang$^{a),b)}$, Jie Xiong$^{c)}$ and Lihu Xu$^{c)}$  }\\
  \footnotesize{$^{a)}$School of Mathematical Sciences, Beijing Normal University, Beijing 100875, China}\\
    \footnotesize{$^{b)}$ Department of Mathematics, Swansea University, Singleton Park, SA2 8PP, UK}\\
   \footnotesize{$^{c)}$ Department of Mathematics,   University of  Macau, Taipa, Macau, China}\\
\footnotesize{Email: \tttext{wangfy@bnu.edu.cn}; \tttext{F.Y.Wang@swansea.ac.uk}; \tttext{jiexiong@umac.mo}; \tttext{lihuxu@umac.mo} }}

\begin{document}
\def\tttext#1{{\normalfont\ttfamily#1}} \def\div{{\rm div}}
\def\R{\mathbb R}  \def\ff{\frac} \def\ss{\sqrt} \def\B{\mathbf B}
\def\N{\mathbb N} \def\kk{\kappa} \def\m{{\bf m}}
\def\dd{\delta} \def\DD{\Delta} \def\vv{\varepsilon} \def\rr{\rho}
\def\<{\langle} \def\>{\rangle} \def\GG{\Gamma} \def\gam{\gamma}
  \def\nn{\nabla} \def\pp{\partial} \def\EE{\scr E}
\def\d{\text{\rm{d}}} \def\bb{\beta} \def\aa{\alpha} \def\D{\scr D}
  \def\si{\sigma} \def\ess{\text{\rm{ess}}}
\def\beg{\begin} \def\beq{\begin{equation}}  \def\F{\scr F}
\def\Ric{\text{\rm{Ric}}} \def\Hess{\text{\rm{Hess}}}
\def\e{\text{\rm{e}}} \def\ua{\underline a} \def\OO{\Omega}  \def\oo{\omega}
 \def\tt{\tilde} \def\Ric{\text{\rm{Ric}}}
\def\cut{\text{\rm{cut}}} \def\P{\mathbb P}
\def\C{\scr C}     \def\E{\mathbb E}
\def\Z{\mathbb Z} \def\II{\mathbb I}
  \def\Q{\mathbb Q}  \def\LL{\Lambda}\def\L{\scr L}
  \def\B{\scr B}    \def\ll{\lambda}
\def\vp{\varphi}\def\H{\mathbb H}\def\ee{\mathbf e}
\def\ep{\scr R}

\maketitle
\begin{abstract} By using the dimension-free Harnack inequality and the integration by parts formula for the associated diffusion semigroup, we prove the central limit theorem, the moderate deviation principle, and the logarithmic iteration law for the sample entropy production rate of  stochastic differential equations with Lipschitz continuous and dissipative drifts.
\end{abstract} \noindent

 \vskip 2cm

\maketitle

\vspace{4mm}

\section{Introduction}
The  entropy production rate (EPR in short) is a key element of the second law of thermodynamics for open systems, see for instance \cite{GM,JNPS14, CGXX, Qian01,TM012}. In this paper we characterize the asymptotic behaviors of the sample EPR for diffusion processes.

Let $(X_t)_{t\ge 0}$ be a stationary diffusion process on $\R^d$ with invariant probability measure $\mu$. It is called reversible if $X_{[0,t]}:=
(X_r)_{0\le r\le t}$ and the reverse $\bar X_{[0,t]}:=(X_{t-r})_{0\le r\le t}$ are identified in distributions for all $t> 0.$ In case that the process is not reversible, the sample EPR   is an important object to  measure the difference between the distributions  of  the process and its reverse. More precisely, for $\P_{[0,t]}$ and $\bar \P_{[0,t]}$ being the distributions of $X_{[0,t]}$ and $\bar X_{[0,t]}$ respectively, the  sample EPR of the process is defined as (see \cite{JQQ})
\
$$
\ep_t(X_{[0,t]})=\ff 1 t \log \ff{\d \P_{[0,t]}}{\d \bar \P_{[0,t]}},\ \ t> 0,
$$
which is a measurable function on $C([0,t];\R^d)$ for every $t>0$. If $\P_{[0,t]}$ is not absolutely continuous with respect to $\bar \P_{[0,t]}$, we set $\ep_t(X_{[0,t]})=\infty.$
It is well known that  $\ep:= \E \ep_t(X_{[0,t]}) $ is non-negative and independent of $t$, and when it is finite we have
\
\beq\label{LM} \lim_{t\to\infty} \ep_t(X_{[0,t]})=\ep\ \ {\rm a.s.}
\end{equation}
according to the ergodic theorem.
The purpose of this paper is to investigate long time behaviors of $\ep_t(X_{[0,t]})$, which include the central limit theorem (CLT in short), the moderate deviation principle (MDP in short) and the logarithmic iteration law (LIL in short).

Consider the following stochastic differential equation (SDE in short) on $\R^d$:
\beq\label{SDE} \d X_t =B(X_t)\d t+ \si \d W_t,\end{equation} where $W_t$ is the $d$-dimensional Brownian motion on a complete filtration probability space $(\OO,\F, \{\F_t\}_{t\ge 0}, \P)$,  $\si$ is an invertible $d\times d$-matrix, and $B:\R^d\to \R^d$ is Lipschitz continuous so that $\nn B$ exists with $\|\nn B\|_\infty<\infty$. We further assume that $B$ satisfies the dissipativity condition
\beq\label{C1} \<B(x)-B(y), x-y\>\le \kk|x-y|-K |x-y|^2,\ \ x,y\in\R^d\end{equation}   for some constants $\kk\ge 0,K>0.$  Note that \eqref{C1} holds   for $B:=B_0+B_1$ where
$B_0$ is bounded   and $B_1\in C^1$ such that $\<\nn_v B_1(x),v\>\le -K |v|^2$ for $x,v\in \R^d$.
It is well known that in this situation   the SDE \eqref{SDE} has a unique non-explosive solution for any initial distributions, and the associate Markov semigroup $P_t$ has a unique invariant probability measure $\mu$. According to \cite{BR}, we have $\mu(\d x)=\rr(x)\d x$ for some strictly positive density function $\rr\in \cap_{p>1}W_{loc}^{p,1}(\d x)$, see Proposition \ref{P2.1} below for details. Throughout the paper, we   denote
$\nu(f)=\int_{\R^d}f\d\nu$ for a measure $\nu$ and $f\in L^1(\nu).$

We now formulate the sample EPR for the solution to \eqref{SDE}.
It is well known that the reverse process is a weak solution to the SDE (see e.g. \cite[Theorem 3.3.5]{JQQ})
\beq\label{RS}\d \bar X_t=   \{\si\si^*\nn\log \rr(\bar X_t)- B(\bar X_t)\}\d t + \si\d W_t.\end{equation}
Since the drift is in $L_{loc}^p(\d x)$ for all $p>1$,   according to   \cite{Zhang}   this SDE  has a unique solution for any initial point.
We will prove
\beq\label{EXP1} \mu\left(\exp[\vv(|B|^2+|\nn\log\rr|^2)]\right)<\infty  {\rm \ \ for \ some \ constant} \ \vv>0\end{equation}
and that   the process
\beq\label{MT}  M_t:=\exp\bigg[\int_0^t \<\si^*\nn\log \rr-2 \si^{-1}B, \d W_s\>
 - \ff 1 2 \int_0^t |\si^*\nn\log \rr-2 \si^{-1}B|^2\d s\bigg],\ \ t\ge 0  \end{equation}   is a martingale (see Proposition \ref{P2.1} below).
 Then by the Girsanov theorem,
$$\bar W_s:= W_s+ \int_0^s\big\{2 \si^{-1}B(X_u)-\si^*\nn\log \rr(X_u)\big\}\d u,\ \ s\in [0,t]$$ is a $d$-dimensional Brownian motion under the probability $\d\Q_t:=M_t\d\P$. Reformulating \eqref{SDE} as
$$\d X_t= \{\si\si^*\nn\log \rr(X_t)- B(X_t)\}\d t + \si\d \bar W_t,$$
we see that
the solution to \eqref{RS} is non-explosive  and,
by the weak uniqueness, we obtain   $$\E f(\bar X_{[0,t]})= \E_{\Q_t} f( X_{[0,t]}) = \E\big\{M_t f(X_{[0,t]})\},\ \   F\in \B_b(C([0,t];\R^d)).$$
This implies  $\ff{\d \bar \P_{[0,t]}}{\d \P_{[0,t]}}(X_{[0,t]})= M_t,$ so that
  the sample EPR of $X_t$ can be  formulated as
\beq\label{EPR}  \beg{split} &\ep_t(X_{[0,t]})=\ff 1 t \log \frac1{M_t}\\
&= -\ff 1 t \int_0^t \<\si^*\nn\log \rr-2 \si^{-1}B, \d W_s\>
 + \ff 1 {2t} \int_0^t |\si^*\nn\log \rr-2 \si^{-1}B|^2\d s.\end{split}\end{equation}

Let $ \scr P(\R^d)$ be  the space of probability measures on $\R^d$. For any $\nu\in \scr P(\R^d)$,   let
 $(X_t^\nu)_{t\ge 0}$ be the solution to \eqref{SDE} with initial distribution $\nu$.  When $\nu=\dd_x$, the Dirac measure at point $x$, we simply denote
 $X^\nu$ by $X^x$. Let
 $$U_\mu^p(l)= \bigg\{\nu\in \scr P(\R^d):\ \int_{\R^d} \Big(\ff{\d\nu}{\d\mu}\Big)^p\d\mu\le l\bigg\},\ \ p>1, l>0.$$
The main result of the paper is the following,  which includes  CLT, MDP and LIL  for
the sample EPR process $\ep_t(X_{[0,t]}).$

\beg{thm}\label{T1.1} Assume  that $B$ is Lipschitz continuous and $\eqref{C1}$  holds for some constants  $\kappa \ge 0$ and $K>0$. Then   the following assertions hold:
\beg{enumerate} \item[$(1)$]  $\eqref{EXP1}$ holds, and  $\dd:= \lim_{t\to\infty} t\E \{\ep_t(X_{[0,t]}^\mu)-\ep\}^2<\infty$ exists.
\item[$(2)$] {\bf (CLT)} For any $p>1$ and  $l>0$, $\lim_{t\to\infty}\P\big(\ss t\{\ep_t(X_{[0,t]}^\nu)-\ep\}\in \cdot\big)=N(0,\dd)$ weakly and uniformly in $\nu\in U_\mu^p(l)$, where $N(0,\dd)$ is the centered Gaussian distribution with variance $\dd$.
\item[$(3)$] {\bf (MDP)} For  any    $\ll: (0,\infty)\to (0,\infty)$ with $\ll(t)\land \ff{\ss t}{\ll(t)}\to\infty$ as $t\to\infty$,   any
measurable set $A\subset \R$ and constants $p>1, l>0$,
 \beg{align*} -\inf_{u\in A^{\rm o} }\ff {u^2}{2\dd} &\le \liminf_{t\to\infty}\ff 1 {\ll(t)^2} \log \P  \bigg(\ff{\ss t}{\ll(t)}\big\{\ep_t(X_{[0,t]}^\nu)-\ep\big\}\in A\bigg)\\
 &\le \limsup_{t\to\infty} \ff1 {\ll(t)^2} \log \P  \bigg(\ff{\ss t}{\ll(t)}\big\{\ep_t(X^\nu_{[0,t]})-\ep\big\}\in A\bigg)\le -\inf_{u\in \bar A } \ff {u^2}{2\dd}\end{align*} holds   uniformly in $\nu\in U_\mu^p(l)$, where $A^{\rm o}$ and $\bar A$ are the interior and closure of $A$ respectively.
 \item[$(4)$] {\bf (LIL)}  For any $\nu\in\scr P(\R^d)$ with $\ff{\d\nu}{\d\mu}\in L^p(\mu)$ for some $p>1$,   $\P$-a.s.
\beg{align*}&\limsup_{t\to\infty} \ff{\ss t}{\ss{2\log\log t}}\big\{\ep_t(X_{[0,t]}^\nu)-\ep\big\}= \ss{\dd},\\
   &\liminf_{t\to\infty} \ff{\ss t}{\ss{2\log\log t}}\big\{\ep_t(X_{[0,t]}^\nu)-\ep\big\}=-\ss\dd.
 \end{align*}
     \end{enumerate}\end{thm}

Since the rate function in Theorem \ref{T1.1}(3) is continuous,    for any domain $A\subset \R$ we have
$$\lim_{t\to\infty} \ff1 {\ll(t)^2} \log \P  \bigg(\ff{\ss t}{\ll(t)}\big\{\ep_t(X^\nu_{[0,t]})-\ep\big\}\in A\bigg)= -\inf_{u\in   A } \ff {u^2}{2\dd}.$$
The next result extends Theorem \ref{T1.1}(2)-(4) to $\nu=\dd_x$, the Dirac measure at point $x$, which is singular with respect to $\mu$.

\beg{thm}\label{C1.2} In the situation of Theorem $\ref{T1.1},$   the following assertions hold.
  \beg{enumerate}
\item[$(1)$]   $\lim_{t\to\infty}\P\big(\ss t\{\ep_t(X_{[0,t]}^x)-\ep\}\in \cdot\big)=N(0,\dd)$ weakly and locally uniformly  in $x\in\R^d.$
\item[$(2)$]   For  any    $\ll: (0,\infty)\to (0,\infty)$ with $\ll(t)\land \ff{\ss t}{\ll(t)}\to\infty$ as $t\to\infty$,   and any measurable set
  $A\subset \R,$
 \beg{align*} -\inf_{u\in A^{\rm o}}\ff {u^2}{2\dd} &\le \liminf_{t\to\infty}\ff 1 {\ll(t)^2} \log \P  \bigg(\ff{\ss t}{\ll(t)}\big\{\ep_t(X_{[0,t]}^x)-\ep\big\}\in A\bigg)\\
 &\le \limsup_{t\to\infty}\ff1 {\ll(t)^2} \log \P  \bigg(\ff{\ss t}{\ll(t)}\big\{\ep_t(X^x_{[0,t]})-\ep\big\}\in A\bigg)\le -\inf_{u\in \bar A } \ff {u^2}{2\dd}\end{align*} holds locally  uniformly  in $x\in\R^d.$
   \item[$(3)$]  For any $x\in\R^d,$    $\P$-a.s.
\beg{align*}&\limsup_{t\to\infty} \ff{\ss t}{\ss{2\log\log t}}\big\{\ep_t(X_{[0,t]}^x)-\ep\big\}= \ss{\dd},\\
   &\liminf_{t\to\infty} \ff{\ss t}{\ss{2\log\log t}}\big\{\ep_t(X_{[0,t]}^x)-\ep\big\}=-\ss\dd.
 \end{align*}
  \end{enumerate}
 \end{thm}

We will prove the above two results   in Section 3,  for which       some  preparations    are presented in  Sections 2. Finally, SDEs with multiplicative noise are  discussed in Section 4.

 \section{Preparations}

Let $P_t$ be the Markov semigroup associated to the SDE \eqref{SDE}, and let
$$\scr U:=\big\{\rr_\nu\d\mu:  \rr_\nu\ge 0, \mu(\rr_\nu)=1, \mu(\rr_\nu^p)<\infty\  \text{for \ some\ }p>1\big\}.$$
 The main result of this section is the following.

 \begin{prp} \label{P2.1}   Assume  that $B$ is Lipschitz continuous and $\eqref{C1}$  holds for some constants  $\kappa \ge 0$ and $K>0$. Then:
 \beg{enumerate} \item[$(1)$] $P_t$ has a unique invariant probability measure $\mu$, which has strictly positive density $\rr\in \cap_{p>1}W^{p,1}_{loc}(\d x),$ and
 $\mu(\e^{\vv (|\cdot|^2+|\nn\log\rr|^2)})<\infty$ holds for some constant $\vv>0.$
  \item[$(2)$] The density   $p_t(x,y)$  of $P_t$ with respect to $\mu$ satisfies
\beq\label{HEAT} \mu(p_t(x,\cdot)^q)\le  \exp\bigg[ 2q(q-1)\kk^2 \|\si^{-1}\|^2  + \ff{4q(q-1)K\|\si^{-1}\|^2(\mu(|\cdot|^2)+|x|^2)}{\e^{2Kt}-1}\bigg]\end{equation} for all $q>1, t>0, x\in\R^d.$
   Consequently, $P_t$ is hyperbounded, i.e. $\|P_t\|_{L^2(\mu)\to L^4(\mu)}<\infty$ for some $t>0.$
  \item[$(3)$]  For any $t>0$, $1$ is a simple eigenvalue of $P_t$.
  \item[$(4)$] For any $p>1\lor\ff d 2$ there exist a constant $c>0$ and  a positive  function $H\in C(\R^d)$ such that
  $$\int_0^t P_s |f|(x)\d s \le \mu(|f|^p)^{\ff 1 p} H(x)\left(t+ t^{\ff{2p-d}{2p}}\right) ,\ \ x\in \R^d, t\ge 0, f\in L^p(\mu).$$
  \item[$(5)$] If $\psi: \R^d\to \R^d$ is measurable such that $\mu(\e^{\vv |\psi|^2})<\infty$ for some $\vv>0$, then for any $\nu\in \{\dd_x: x\in\R^d\}\cup\scr U$,
  $$M_t^\nu := \exp\bigg[\int_0^t \<\psi(X_s^\nu), \d W_s\> -\ff 1 2 \int_0^t |\psi(X_s^\nu)|^2\d s\bigg],\ \ t\ge 0$$ is a martingale.
 \end{enumerate} \end{prp}

 To prove this result we need the following lemma on exponential integrability,   integration by parts formula and Harnack inequality.

 \beg{lem} \label{L2.1}   Assume that $B$ is Lipschitz continuous and  $\eqref{C1}$  holds for some constants $\kk\ge 0$ and $K>0$. Then:
 \beg{enumerate} \item[$(1)$] There exist constants $\vv, c >0$ such that
 \beq\label{EST} \E \int_0^t \e^{\vv|X_s^x|^2} \d s \le c\big(t+  \e^{\vv |x|^2}\big),\ \ x\in\R^d, t\ge 0.\end{equation}
 \item[$(2)$] For any $f\in \B_b^+(\R^d)$, $T>0, p>1$ and $x,y\in\R^d$,
 \beq\label{HI} (P_Tf(x))^p\le (P_Tf^p(y))\exp\bigg[\ff{2p\kk^2\|\si^{-1}\|^2(\e^{KT}-1)}{(p-1)(\e^{KT}+1)}+  \ff{2pK\|\si^{-1}\|^2|x-y|^2}{(p-1)(\e^{2KT}-1)}\bigg].\end{equation}
 \item[$(3)$]  For any   $f\in C^1_b(\R^d)$,
 \beq\label{ITF} P_T \nn f(x)= \E \bigg[\ff{f(X_T^x)}T \int_0^T \big\{\si^{-1}(I-t \nn B(X_t^x))\big\}^* \d W_t\bigg],\ \ T>0, x\in\R^d. \end{equation}\end{enumerate}
   \end{lem}
 \beg{proof} Assertion (1) follows from condition \eqref{C1} by It\^o's formula, and the other two can be easily proved by using coupling by change of measures as in  \cite{Wan13}.
 We include below   brief proofs of these assertions for completeness.

(1) It is easy to see that the generator of the diffusion process is
$$L =\ff 1 2 \sum_{i,j=1}^d(\si\si^*)_{ij}\pp_i\pp_j + \sum_{i=1}^d B_i(x)\pp_i.$$ Then condition \eqref{C1} implies that for small enough $\vv>0,$
\begin{equation}\label{2.5}
L\e^{\vv |\cdot|^2}(x) \le   C_1-C_2  \e^{\vv |x|^2},\ \ x\in\R^d
\end{equation}
holds for some constants $C_1, C_2>0$.
So, by It\^o's formula we obtain
$$C_2 \E\int_0^t \e^{\vv|X_s^x|^2} \d s \le \e^{\vv |x|^2} + C_1 t,$$ which implies  \eqref{EST} for $c:= \ff{1\lor C_1}{C_2}.$

(2) For fixed $x,y\in\R^d$ and $T>0$, let $X_t=X_t^x$ solve \eqref{SDE} for $X_0=x$, and  construct $Y_t$ with $Y_0=y$ as follows. For
$$\xi_t:= \kk \e^{-K(T-t)}+ \ff{2K\e^{Kt}|x-y|}{\e^{2KT}-1},\ \ t\ge 0,$$     the SDE
$$\d Y_t= \bigg\{B(Y_t)+ \xi_t \ff{X_t-Y_t}{|X_t-Y_t|}\bigg\}\d t +\si\d W_t,\ \ Y_0=y$$ has a unique solution before  the coupling time
$$\tau:=\inf\{t\ge 0:\ X_t=Y_t\}.$$ Take $Y_t=X_t$ for $t\ge \tau$. Then $(Y_t)_{t\ge 0}$ solves the SDE
\beq\label{S'} \d Y_t= B(Y_t)\d t    +\si\d \tt W_t,\ \ Y_0=y\end{equation} for
$$\tt W_t:= W_t + \int_0^{t\land \tau}  \xi_s\ff{\si^{-1}(X_s-Y_s)}{|X_s-Y_s|}\d s,\ \ s\ge 0.$$ By the Girsanov theorem,
$(\tt W_t)_{t\in [0,T]}$ is a $d$-dimensional Brwonian motion under the probability $\Q:= R\P$ for
$$R:= \exp\bigg[-\int_0^{T\land\tau} \Big\< \xi_s\ff{\si^{-1}(X_s-Y_s)}{|X_s-Y_s|}, \d W_s\Big\>- \ff 1 2\int_0^{T\land \tau} \Big|  \xi_s\ff{\si^{-1}(X_s-Y_s)}{|X_s-Y_s|}\Big|^2\d s\bigg].$$
Combining \eqref{S'} with \eqref{SDE} and using condition \eqref{C1}, we obtain
$$\d|X_t-Y_t|\le (\kk- \xi_t-K |X_t-Y_t|)\d t,\ \ t\le \tau\land T.$$ This together with the definition of $\xi_t$ leads to
\beg{align*}|X_t-Y_t|&\le\ \e^{-Kt}|x-y| + \int_0^t \e^{-K(t-s)}(\kk-\xi_s)\d s\\
&= \ff{\kk\e^{-Kt}}{K}\Big(\e^{Kt}-1-\ff{\e^{2Kt}-1}{\e^{KT}+1}\Big)+ |x-y|\e^{-Kt} \Big( 1- \ff{\e^{2Kt}-1}{\e^{2KT}-1}\Big),\ \ t\in [0, T\land\tau].\end{align*} So, if $\tau>T$ then  by the definition of $\tau$ we have
$$0<|X_T-Y_T|\le \ff{\kk\e^{-KT}}{K}\Big(\e^{KT}-1-\ff{\e^{2KT}-1}{\e^{KT}+1}\Big)+ |x-y|\e^{-KT} \Big( 1- \ff{\e^{2KT}-1}{\e^{2KT}-1}\Big)=0,$$ which is impossible. Therefore, $\tau\le T$ a.s. so that $X_T=Y_T.$ Combining this with \eqref{S'} and noting that $\tt W_t$ is  a Brownian motion under $R\P$, for $f\in \B_b^+(\R^d)$ we have
\
\begin{equation*}
(P_Tf(y))^p  = (\E \{R f(Y_T)\})^p = (\E\{Rf(X_T)\})^p \le \{\E f^p(X_T)\}(\E R^{\ff p{p-1}})^{p-1}
\end{equation*}
and
\
\beg{align*}
\E R^{\ff p{p-1}} &\le \E\exp\bigg[\ff{p\|\si^{-1}\|^2}{2(p-1)^2}\int_0^T |\xi_s|^2\d s  - \ff{p^2}{2(p-1)^2}
\int_0^{T\land \tau}  \Big|\xi_s\ff{\si^{-1}(X_s-Y_s)}{|X_s-Y_s|} \Big|^2\d s\\
&\qquad\qquad \qquad\qquad \qquad\qquad \qquad -\ff p{p-1}\int_0^{T\land\tau}  \Big\< \xi_s\ff{\si^{-1}(X_s-Y_s)}{|X_s-Y_s|}, \d W_s \Big\>\bigg]\\
&=\exp\bigg[\ff{p\|\si^{-1}\|^2}{2(p-1)^2}\int_0^T |\xi_s|^2\d s\bigg]\\
&\le \exp\bigg[\ff{2p\kk^2\|\si^{-1}\|^2(\e^{KT}-1)}{(p-1)^2(\e^{KT}+1)}+ \ff{2pK\|\si^{-1}\|^2|x-y|^2}{(p-1)^2(\e^{2KT}-1)}\bigg].
\end{align*}
Hence,
\
$$(P_Tf(y))^p \le  (P_T f^p(x)) \exp\bigg[\ff{2p\kk^2\|\si^{-1}\|^2(\e^{KT}-1)}{(p-1)(\e^{KT}+1)}+ \ff{2pK\|\si^{-1}\|^2|x-y|^2}{(p-1)(\e^{2KT}-1)}\bigg].$$

(3) Again let $X_t=X_t^x$. For any $v\in\R^d$ with $|v|=1$ and $r\in [0,1]$, let $Y_t^r$ solve the SDE
$$\d Y_t^r = \Big\{B(X_t)+ \ff {r v}T\Big\}\d t + \si\d W_t,\ \ Y_0^r=x.$$Then
\beq\label{CP} Y_t^r-X_t= \ff{r t v  }T,\ \ t\in[0,T].\end{equation}
Since $B$ is Lipschitz continuous, we have
$$|B(Y_t^r)-B(X_t)|\le\|\nn B\|_\infty |Y_t^r-X_t|\le \|\nn B\|_\infty r|v|<\infty,\ \ t\in [0,T], r\in [0,1].$$
Then by the Girsanov theorem,
$$W_t^r:= W_t+ \int_0^t  \si^{-1}\Big\{B(X_s)-B(Y_s^r) + \ff{rv}T\Big\} \d s,\ \ s\in [0,T]$$ is a $d$-dimensional Brownian motion under the probability $R_r\P$, where
$$R_r:= \exp\bigg[-\int_0^T \Big\< \si^{-1}\Big\{B(X_s)-B(Y_s^r) + \ff{rv}T\Big\},\d W_s\Big\>-\ff 1 2 \int_0^T\Big|\si^{-1}\Big\{B(X_s)-B(Y_s^r) + \ff{rv}T\Big\} \Big|^2\d s\bigg].$$
Combining this with  $Y_T^r=X_T+ rv$ due to \eqref{CP}, we obtain
$$P_t f(x)= \E\{R_r f(Y_T^r)\}= \E\{R_r f(X_T+rv)\},\ \ r\in [0,1].$$
Due to   \eqref{CP}, $\|\nn B\|_\infty<\infty$ and the definition of $R_r$, for any $f\in C_b^1(\R^d)$ we may take derivative for both sides in $r$ at $r=0$ to derive
\beg{align*} 0&= \ff{\d}{\d r}\Big|_{r=0} P_T f(x)= \E \Big\{f(X_T)\ff{\d R_r}{\d r}\Big|_{r=0}\Big\}+ \E\{(\nn_v f) (X_T)\}\\
&= \ff 1 T \E\bigg\{f(X_T)  \int_0^T \<\si^{-1}\{t\nn_vB(X_s)- v\},\d W_s\>\bigg\} + P_T(\nn_vf)(x),\ \ v\in\R^d.\end{align*} This implies \eqref{ITF}
  \end{proof}

\beg{proof}[Proof of Proposition \ref{P2.1}] (1)  It is well known that  \eqref{2.5} implies the existence of invariant probability measure and that any invariant  probability measure $\mu$ satisfies    $\mu(\e^{\vv|\cdot|^2})<\infty$.
 By the Harnack inequality \eqref{HI}, $\mu$ is the unique invariant probability measure (see  \cite[Theorem 1.4.1(3)]{Wan13}
or \cite[Proposition 3.1]{WaYu11}). As already indicated in the Introduction that according to \cite{BR}, $\mu(\d x)=\rr(x)\d x$ holds for some strictly positive
$\rr\in \cap_{p>1} W_{loc}^{p,1}(\d x)$. It remains to prove that $\mu(\e^{\vv |\nn\log\rr|^2})<\infty$ for some $\vv>0.$

 Let $X_t^\mu$ be the solution to \eqref{SDE} with initial distribution $\mu$.
Since $\mu$ is $P_t$-invariant, by taking integral for \eqref{ITF} with respect to $\mu(\d x)$ we obtain
\beg{align*} \mu(\nn f) &=\mu(P_1\nn f) =   \E  \bigg\{ f(X_1^\mu)   \int_0^1 \{\si^{-1}(I-t \nn B(X_t^\mu))\}^*\d W_t\bigg\}\\
&=  \E   \bigg\{ f(X_1^\mu)  \E\bigg(\int_0^1 \{I-\si^{-1}(t \nn B(X_t^\mu))\}^*\d W_t\bigg|X_1^\mu\bigg)\bigg\},\ \ f\in C_0^1(\R^d).\end{align*}
On the other hand, by the integration by parts formula for the Lebesgue measure, for any $f\in C_0^1(\R^d)$ we have
$$\mu(\nn f)= - \mu(f\nn\log\rr)= -\E\big\{f(X_1^\mu)\nn\log\rr(X_1^\mu)\big\}.$$ Combining this with the above display we obtain
$$\nn\log\rr(X_1^\mu)= \E\bigg(\int_0^1 [\si^{-1}(t \nn B(X_t^\mu)-I)]^*\d W_t\bigg|X_1^\mu\bigg),\ \  {\rm a.s.}$$
Then by Jensen's inequality and noting that $\|\nn B\|_\infty<\infty$,  we have
\beg{align*} \mu(\e^{\vv |\nn\log \rr|^2})&=\E \exp\big[\vv |\nn \log \rr(X_1^\mu)|^2\big] \\
& \le \E\bigg\{\E\bigg[\exp\bigg(\vv\bigg|\int_0^1 [\si^{-1}(t \nn B(X_t^\mu)-I)]^*\d W_t\bigg|^2\bigg)\bigg|X_1^\mu\bigg]\bigg\}\\
&= \E \exp\bigg[\vv \bigg|\int_0^1 \{\si^{-1}(t \nn B(X_t^\mu)-I)\}^*\d W_t\bigg|^2\bigg]<\infty
\end{align*}
 for small enough $\vv>0$.

(2) By the Harnack inequality \eqref{HI},  for $x,y\in\R^d$ and $t>0$ we have
$$\bigg|\int_{\R^d} p_t(x,z) f(z) \mu(\d z)\bigg|^{\ff q {q-1}} \exp\Big[-\ff{2qK \|\si^{-1}\|^2|x-y|^2}{\e^{2Kt}-1}\Big]\le \exp\big[2q\kk^2\|\si^{-1}\|^2\big] P_t|f|^{\ff q {q-1}}(y).$$
Taking integral with respect to $\mu(\d y)$, when $\mu(|f|^{\ff q {q-1}})=1$ we obtain
$$\bigg|\int_{\R^d} p_t(x,z)f(z)\mu(\d z)\bigg|^{\ff q {q-1}} \le \ff{\exp[2q\kk^2\|\si^{-1}\|^2])}{\mu \big(\exp\big[-\ff {2qK\|\si^{-1}\|^2|x-\cdot|^2}{\e^{2Kt}-1}\big]\big)}.$$
This implies
\beq\label{EY}  \mu(p_t(x,\cdot)^q)\le \bigg\{\ff{\exp[2q\kk^2\|\si^{-1}\|^2]}{\mu \big(\exp\big[-\ff {2qK\|\si^{-1}\|^2|x-\cdot|^2}{\e^{2Kt}-1}\big]\big)}\bigg\}^{q-1}. \end{equation}
By Jensen's inequality we have
\beg{align*}&\mu\bigg(\exp\Big[-\ff {2qK\|\si^{-1}\|^2|x-\cdot|^2}{\e^{2Kt}-1}\Big]\bigg)
 \ge \exp\bigg[- \mu\Big(\ff {2qK\|\si^{-1}\|^2|x-\cdot|^2}{\e^{2Kt}-1}\Big)\bigg]\\
  & \ge\exp\bigg[-\ff{4qK\|\si^{-1}\|^2(|x|^2+\mu(|\cdot|^2)}{\e^{2Kt}-1}\bigg].\end{align*}
Substituting this into \eqref{EY} we prove \eqref{HEAT}.

Next, by \eqref{HEAT}, there exists a constant $C>0$ such that  if $\mu(f^2)\le 1$ then
$$|P_t f(x)|^4 \le [\mu(p_t(x,\cdot)^2)]^2[\mu(f^2)]^2 \le \mu(p_t(x,\cdot)^4) \le  \e^{\ff {C(1+t +|x|^2)} {t}}.$$
 Since  $\mu(\e^{\vv|\cdot|^2})<\infty$ for some constant $\vv>0$, when $t>0$ is large enough this implies
 $$\sup_{\mu(f^2)\le 1}\mu((P_tf)^4)\le \e^{\ff {C(1+t)}{t}}\mu(\e^{\ff{C|\cdot|^2}{t}})<\infty.$$ Thus, $P_t$ is hyperbounded.

 (3) Since due to (1) we have $|\nn\log\rr|+|B|\in L^p(\mu)$ for any $p>1$,  by   \cite[Proposition 2.11]{MR} with $c=b=0$ and $d= B-\si\si^*\nn\log\rr$ which is divergence free as $\mu$ is an invariant probability measure,   $P_t$ is associated to a Dirichlet form    with  symmetric part
  $$ \EE(f,g):= \mu(\<\si\si^*\nn f,\nn g\>), \ \ f,g\in H_\si^{2,1}(\mu),$$ where
  $H_\si^{2,1}(\mu)$ is the completion of $C_0^\infty(\R^d)$ under the  Sobolev norm
  $$\|f\|_{2,1}:=\ss{\mu(f^2)+\EE(f,f)}.$$
  Obviously, the Dirichlet form is irreducible so that $P_t\to\mu$ in $L^2(\mu)$ as $t\to \infty.$ If $f\in L^2(\mu)$ such that $P_t f=f$ for some $t>0$, then $f=P_{nt}f\to \mu(f)$  in $L^2(\mu)$ as  $n\to\infty$,  so that $f$ has to be  constant. Thus,   $1$ is a simple eigenvalue of $P_{t}$.

 (4)   By the Harnack inequality \eqref{HI}   we obtain
 \beq\label{1*}  \big(P_s |f| (x)\big)^{p} \int_{\R^d} \e^{-c-c|x-y|^2/s}\mu(\d y)\le \mu(|f|^p)  \end{equation} for some constant $c=c(p)>0$.
  Since $\rr\in C(\R^d)$ is strictly positive as already explained in Introduction due to \cite{BR}, we have $\mu(B(x, \ss s))
   \ge (1\land s)^{d/2} h_1(x)$ for some strictly positive $h_1\in C(\R^d)$ and all $s\ge 0, x\in\R^d$. Then
 $$\int_{\R^d} \e^{-c-c|x-y|^2/s}\mu(\d y)\ge \e^{-c} \mu\big(B(x,\ss s)\big) \ge \e^{-c} (1\land s)^{d/2}h_1(x),\ \ s> 0, x\in \R^d.$$ Combining this with \eqref{1*}
 we obtain
 $$P_s  |f| (x)\le \mu(|f|^p)^{\ff 1 p}   (1\land s)^{-\ff d{2p}} h_2(x),\ \ s>0, x\in \R^d$$ for some positive $h_2\in C(\R^d).$ Since $p>\ff d 2,$  this implies   the desired estimate.

   (5) By (4) and $\mu(\e^{\vv|\psi|^2})<\infty$ we have
   $\E\int_0^t |\psi(X_s^x)|^2\d s<\infty$ for any $t>0$ and $x\in \R^d$, and for $\rr_\nu:= \ff{\d\nu}{\d\mu}\in L^q(\mu)$,
 \begin{equation*}
 \begin{split}
 \E\int_0^t |\psi(X_s^\nu)|^2\d s&=  \int_0^t \mu(\rr_\nu P_s |\psi|^2) \\
  &\le \mu(\rr_\nu^q)^{\ff 1 q} \int_0^t \left[\mu(P_s |\psi|^{\ff{2q}{q-1} })\right]^{1-\frac 1q} \d s \\
 &=\mu(\rr_\nu^q)^{\ff 1 q} \left[\mu( |\psi|^{\ff{2q}{q-1} })\right]^{1-\frac 1q} <\infty.
\end{split}
 \end{equation*}
 Then for $\nu= \dd_x$ or $\nu\in \scr U$, $M_t^\nu$ is a well defined supermartingale.  It suffices to prove $\E M_t^\nu=1$ for any $t\ge 0$.
   Since $\E M_t^\nu = \int_{\R^d}(\E M_t^{x})\nu(\d x)$, it remains to show that $(M_t^x )_{t\ge 0}$ is a martingale for any $x\in\R^d$. By the Markov property and the Girsanov theorem, this follows if we can find a constant $t_0>0$ such that  the Novikov condition
  \
 \begin{equation}  \label{e:NovCon}
   \E\exp\left(\ff 1 2 \int_0^{t_0} |\psi(X_s^x)|^2\d s\right)<\infty, \ \ \ \forall \ x \in \R^d.
  \end{equation}
 Indeed,   \eqref{e:NovCon} implies that $(M^x_t)_{0 \le t \le t_0}$ is a martingale for all $x \in \R^d$, so that by the Markov property,    for any $s\ge 0$ the process
 $(M_{s,t}^x)_{t\in [s,s+t_0]}:= (\ff{M_t^x}{M_s^x})_{t\in [s,s+t_0]}$ is a martingale under the conditional probability given $X_s^x$. Thus, by induction and the Markov property we prove that $(M_t^x)_{t\ge 0}$ is a martingale for all $x\in \R^d$ as follows: if $(M^x_t)_{0 \le t \le nt_0}$ is a martingale  for some $n\ge 1,$
  then for any $nt_0\le s<t \le (n+1)t_0$ we have
 $$ \E (M_t^x|\F_s)= M_s^x \E(M_{s,t}^x|\F_s)=M_s^x \E(M_{s,t}^x|X_s^x)   =M_s^x.$$

To prove \eqref{e:NovCon},   we take  $t_0=\frac{\vv}d.$  By taking   $p=2d$ in Proposition \ref{P2.1}(4) and using    Jensen's inequality we obtain
   \
 \begin{equation*}
\begin{split}
& \E\exp\left(\ff 1 2 \int_0^{t_0} |\psi(X_s^x)|^2\d s\right)  \le \ff 1 {t_0} \int_0^{t_0} P_s \e^{\ff{t_0} 2 |\psi|^2}(x)\d s \\
 & \le C(t_0)H(x)\left[\mu(\e^{d t_0 |\psi|^2})\right]^{\ff 1 {2d}} <\infty
 \end{split}
  \end{equation*}
  for some constant $C(t_0)>0.$

   \end{proof}

 \section{Proofs of Theorems \ref{T1.1} and \ref{C1.2}}

We will prove the following more general result Theorem \ref{T3.1},  which implies Theorem \ref{T1.1} for $S_t^\nu:= t \big\{\scr R_t(X_{[0,t]}^\nu)-\scr R\big\}$ according to Proposition \ref{P2.1}.

In general, let $X_t$ be a time-homogenous continuous Markov process on $\R^d$ with respect to the filtration $\F_t$ such that the associated Markov semigroup $P_t$ has a unique invariant probability measure $\mu$.
Let  $\psi: \R^d\to\R^d$ be measurable such that $\mu(|\psi|^p)<\infty$ for any $p>1$. Since $\mu$ is $P_t$-invariant,   for any
$\nu\in \scr U$ and any $q>1$, the process $X_s^\nu$ starting at distribution $\nu$ satisfies
\beq \label{U00}
\beg{split}\int_r^t \E |\psi(X_s^\nu)|^q\d s&= \int_r^t \nu(P_s|\psi|^{q})\d s \\
&\le \int_r^t \mu(\rr_\nu^p)^{\ff 1 p} \mu\big(P_s|\psi|^{\ff{2pq}{p-1}}\big)^{1-\frac 1p}\d s\\
& =(t -r)\mu(\rr_\nu^p)^{\ff 1 p} \mu\big(|\psi|^{\ff{2pq}{p-1}}\big)^{1-\frac 1p}
<\infty,\ \ \ t\ge r\ge 0. \end{split}\end{equation}  In particular,  the additive functional
$$S_t^\nu:= \int_0^t\<\psi(X_s^\nu),\d W_s\> +\int_0^t \big\{|\psi(X_s^\nu)|^2-\mu(|\psi|^2)\big\}\d s,\ \ t\ge 0, \nu\in \scr U$$ is well defined.

\beg{thm} \label{T3.1} In the above framework, let $\psi: \R^d\to\R^d$ be measurable such that $\mu(\e^{\vv|\psi|^2})<\infty$ for some constant $\vv>0$.   If $P_t$ is hyperbounded and
there exists $t>0$ such that $1$ is a simple eigenvalue of $P_t$ in $L^2(\mu)$, then the following assertions hold:
\beg{enumerate} \item[$(1)$]  $\dd:= \lim_{t\to\infty}  \ff 1 t \E |S_t^\mu|^2<\infty$ exists.
\item[$(2)$]   For any $p>1$ and  $l>0$, $\lim_{t\to\infty}\P\Big(\ff{S_t^\nu}{\ss t}\in \cdot\Big)=N(0,\dd)$ weakly and uniformly in $\nu\in U_\mu^p(l)$.
\item[$(3)$]   For  any    $\ll: (0,\infty)\to (0,\infty)$ with $\ll(t)\land \ff{\ss t}{\ll(t)}\to\infty$ as $t\to\infty$,    any measurable set
  $A\subset \R$ and constants $p>1, l>0$,
 \beg{align*} -\inf_{x\in A^{\rm o}}\ff {x^2}{2\dd} &\le \liminf_{t\to\infty} \ff 1 {\ll(t)^2} \log \P  \bigg(\ff{ S_t^\nu}{\ll(t) \ss{t}}  \in A\bigg)\\
 &\le \limsup_{t\to\infty} \ff1 {\ll(t)^2} \log \P \bigg(\ff{ S_t^\nu}{\ll(t) \ss{t}}  \in A\bigg)\le -\inf_{x\in \bar A } \ff {x^2}{2\dd}\end{align*} holds   uniformly in $\nu\in U_\mu^p(l)$.
  \item[$(4)$]  If $\mu(|\cdot|^p)<\infty$ for any $p>1$, then   for any $\nu\in\scr P(\R^d)$ with $\rr_\nu:=\ff{\d\nu}{\d\mu}\in L^q(\mu)$ for some $q>1$,   $\P$-a.s.
$$  \limsup_{t\to\infty} \ff{S_t^\nu}{\ss{2t\log\log t}}  \le\ss{\dd},\ \   \liminf_{t\to\infty} \ff{S_t^\nu}{\ss{2t\log\log t}} \ge-\ss\dd.$$
If moreover  $P_t$ has density $p_t(x,\cdot)$ with respect to $\mu$ such that
  \beq\label{CV} \sup_{l\ge 1} \sup_{|x|\le l^{r_0}}\mu(p_l(x,\cdot)^q)<\infty\end{equation} holds for some  $q>1$ and $r_0>0$,    then the equalities hold.

\end{enumerate}\end{thm}

\beg{proof}
(a) According to  \cite[Theorem $2.4'$]{Wu}, the assertions (1)-(3)  hold  provided there exist two constants $t_0,\vv>0$ such that $\sup_{t\in [0,t_0]} \E \e^{\vv |S_t^\mu|}<\infty.$
 For any $r>0$ we have
\beq\label{Q1} \E \e^{r |S_t^\mu|}\le \E \e^{r S_t^\mu }+ \E \e^{-r S_t^\mu }.\end{equation}
Since $\mu(\e^{\vv|\psi|^2})<\infty$, by Schwartz's and Jensen's inequalities,  and noting that $\mu$ is the invariant probability measure, for $t_0:=\ff\vv 4 >0$ we have
\
\beg{align*}
\left(\E \e^{S_t^\mu }\right)^2 &\le \E \exp\left(2\int_0^t\<\psi(X_s^\mu),\d W_s\>-2 \int_0^t |\psi(X_s^\mu)|^2\d s\right)  \\
& \ \ \ \times
 \E\exp\left(4\int_0^t |\psi(X_s^\mu)|^2\d s- 2t\mu(|\psi|^2)\right) \\
&\le \ff 1 t \int_0^t \E\exp\left(4 t |\psi(X_s^\mu)|^2\right) \d s \\
& \le \mu(\e^{4 t_0  |\psi|^2})<\infty,\ \ t\in [0,t_0].
\end{align*}
The same estimate holds for $\E \e^{- S_t^\mu }$.   Then $\sup_{t\in [0,t_0]} \E \e^{ |S_t^\mu|}<\infty$  for some constant $t_0>0.$

(b) To prove   (4), we need the following assertion:   for  any $\nu\in \scr U$ and $p\ge 2$,  there exists a constant $c >0$ such that
\beq\label{L3.1} \E|S_t^\nu-S_s^\nu|^p \le c  (t-s)^{\ff p 2 } \left[\log\log (t-s+\e^2)\right]^{\ff{p}2},\ \ t\ge s\ge 0.\end{equation}
  By \eqref{U00}, it suffices to prove the estimate  for $t-s\ge \e^2$.

 We first consider the case that $\nu=\mu$.  In this case we only need to consider $s=0$ due to the stationary property.
 Since $\mu$ is $P_t$-invariant and $\mu(|\psi|^q)<\infty$ for any $q>1$,    for any $p\ge 2$ there exists a constant $c_1>0$ such that
 \beq\label{U0} \E |S_t^\mu|^p\le c_1 t^{p},\ \ t\ge \e^2.\end{equation}
   Applying Theorem \ref{T1.1}(3) to
 $\ll(t):= \ss{(1+\dd)p\log t}$ for   $t\ge \e,$   we obtain
 $$\P\big(|S_t^\mu|>\ll(t)\ss{2t}\big)\le c_2\e^{-p\log t}= c_2t^{-p},\ \ t\ge \e$$ for some constant $c_2>0.$ Combining this with \eqref{U0}   we arrive at
 \beq\label{U1}\beg{split}  \E|S_t^\mu|^p&\le \E\left(|S_t^\mu|^p1_{\{|S_t^\mu|>\ll(t)\ss{2t}\}}\right) + (2t\ll(t)^2)^{\ff p 2}  \\
 &\le \ss{\P\big(|S_t^\mu|>\ll(t)\ss{2t}\big)\E|S_t^\mu|^{2p}} +c_3 t^{\ff p 2}\ll(t)^p\le c_4t^{\ff p 2} (\log t)^{\ff p2},\ \ t\ge \e\end{split}\end{equation}
 for some constants $c_3,c_4>0.$ Moreover, applying Theorem \ref{T1.1}(3) to
 $\ll(t):= \ss{(1+\dd)p\log\log t}$ for   $t\ge \e^2,$   we obtain
 $$\P\big(|S_t^\mu|>\ll(t)\ss{2t}\big)\le c_5\e^{-p\log\log t}= c_5(\log t)^{-p},\ \ t\ge \e^2$$ for some constant $c_5>0.$ Combining this with \eqref{U1} that
 $\E|S_t^\mu|^{2p}\le c t^p (\log t)^p$ for some constant $c>0$ and $t\ge \e$, we arrive at
 \
\beg{align*}  \E|S_t^\mu|^p&\le \E\left(|S_t^\mu|^p1_{\{|S_t^\mu|>\ll(t)\ss{2t}\}}\right) + (2t\ll(t)^2)^{\ff p 2}  \\
 &\le \ss{\P\big(|S_t^\mu|>\ll(t)\ss{2t}\big)\E|S_t^\mu|^{2p}} +c_6 t^{\ff p 2}(\log\log t)^{\ff p2}\le c_7t^{\ff p 2} (\log\log t)^{\ff p2},\ \ t\ge \e^2\end{align*}
 for some constants $c_5, c_6, c_7>0.$  Thus, the assertion holds for $\nu=\mu$.

 Next, let $\nu\in\scr U$ with $\mu(\rr_\nu^q)<\infty$ for some $q>1$.  By the estimate on $\E |S_t^\mu-S_s^\mu|^p$ we have $\E |S_t^x-S_s^x|^p<\infty, \mu$-a.e. $x$, where $S_t^x:=S_t^\nu$ for $\nu=\dd_x$. Moreover,
 \
\beg{align*}
\E |S_t^\nu-S_s^\nu|^p&=\int_{\R^d} \rr_\nu(x) \E |S_t^x-S_s^x|^p\mu(\d x) \\
& \le \left[\mu(\rr_\nu^q)\right]^{\ff 1 q} \left[\mu\left(\E |S_t^x-S_s^x|^p\right)^{\ff q{q-1}}\right]^{\ff{q-1}q}
 \\
 &\le \left[\mu(\rr_\nu^q)\right]^{\ff 1 q}  \left[\mu\left(\E |S_t^x-S_s^x|^{\ff {pq}{q-1}}\right)\right]^{\ff{q-1}q} \\
 &=\left[\mu(\rr_\nu^q)\right]^{\ff 1 q}  \left(\E |S_t^\mu-S_s^\mu|^{\ff {pq}{q-1}}\right)^{\ff{q-1}q}\\
& \le c (t-s)^{\ff p 2}\{\log\log (t-s+\e^2)\}^{\ff p 2}\end{align*}
 holds for some constant $c>0$.

 (c) To prove   (4),   we will take a sequence $t_n\uparrow \infty$ to replace the continuous limit for $t\uparrow\infty$.   Unlike   the standard choice $t_n=p^n$ for $p>1$ in the literature (see \cite{Chung}),   we   take   $t_n= \e^{n^\theta}$ for some $\theta\in (0,1)$ where  $\theta<1$  is crucial in the argument.

Since $\theta<1$, we may take $p>1$ such that $p(1-\theta)>1.$
For any $\vv>0$, by the stationary property of the process,
   the Burkhold inequality and  \eqref{L3.1}, we obtain
 \beg{align*}
&\P\bigg(\max_{t\in [t_n, t_{n+1}]} |S_t^\nu-S_{t_n}^\nu|> \vv\ss{2t_n\log\log t_n}\bigg)
  \le \ff{c_1\E|S_{t_{n+1}}^\nu-S_{t_n}^\nu|^{2p}}{(2\vv^2t_n\log\log t_n)^p} \\
  & \le \ff{c_2(t_{n+1}-t_n)^p\left[\log\log(t_{n+1}+\e^2)\right]^p}{t_n^p(\log\log t_n)^p}
  \le  c_3n^{(\theta-1)p },\ \ n\ge 2\end{align*}
  for some constants $c_1, c_2,c_3>0$. Hence,
 $$\P\bigg(\max_{t\in [t_n, t_{n+1}]}\ff{|S_t^\nu-S_{t_n}^\nu|}{\ss{2t_n\log\log t_n}}>\vv\bigg) \le c_3n^{(\theta-1)p},\ \ n\ge 2. $$
  Since  $p(1-\theta)>1$, this implies
 $$\sum_{n=2}^\infty \P\bigg(\max_{t\in [t_n, t_{n+1}]}\ff{|S_t^\nu-S_{t_n}^\nu|}{\ss{2t_n\log\log t_n}}>\vv\bigg) <\infty,$$
so that by the Borel-Cantelli lemma,
 $$\limsup_{n\to\infty} \max_{t\in [t_n, t_{n+1}]}\ff{|S_t^\nu-S_{t_n}^\nu|}{\ss{2t_n\log\log t_n}}\le \vv,\ \ \text{a.s.}$$
 By the arbitrariness of $\vv>0$, we arrive at
\beq\label{3*} \limsup_{n\to\infty} \max_{t\in [t_n, t_{n+1}]}\ff{|S_t^\nu-S_{t_n}^\nu|}{\ss{2t_n\log\log t_n}}=0,\ \ \text{a.s.}\end{equation}

(d) We now prove assertion (4) for $\dd=0$. In this case, for any $\vv>0$,  Theorem \ref{T3.1}(3) implies
 $$\lim_{n\to\infty} \ff 1{\log\log t_n}\log \P\Big(\ff{|S_{t_n}^\nu|}{\ss{2t_n\log\log t_n}}>\vv\Big) =-\infty,$$
 so that we may find a constant $c>0$ such that
 $$\P\Big(\ff{|S_{t_n}^\nu|}{\ss{2t_n\log\log t_n}}>\vv\Big)\le c \exp\Big[-\ff 2 \theta \log\log t_n\Big]= \ff c{n^2},\ \ n\ge 1.$$
 Then $\sum_{n=1}^\infty \P\big(\ff{|S_{t_n}^\nu|}{\ss{2t_n\log\log t_n}}>\vv\big)<\infty$. By the Borel-Cantelli lemma, this implies
 $$\limsup_{n\to\infty} \ff{|S_{t_n}^\nu|}{\ss{2t_n\log\log t_n}}\le \vv,\ \ \text{a.s.}$$
 Since $\vv>0$ is arbitrary, we have
 $$\limsup_{n\to\infty} \ff{|S_{t_n}^\nu|}{\ss{2t_n\log\log t_n}}=0,\ \ \text{a.s.}$$
 Combining this with \eqref{3*} we prove (4) for $\dd=0.$

(e)  Let $\dd\in (0,\infty)$. In this case by using $\dd^{-1/2}\psi$ to replace $\psi$, we may and do assume that $\dd=1$. We will only prove the first limit as that of the second is completely similar.  We first prove the upper bound estimate
 \beq\label{U4} \limsup_{t\to\infty} \ff{S_t^\nu}{\ss{2t\log\log t}}  \le 1.\end{equation}
   For any $r\in (0,1)$,  take  $\theta\in ((1+ r)^{-2},1)$ and  $t_n=\e^{n^\theta}$ for $n\ge 1.$
   We have
  $t_{n+1}-t_n\le c_1 t_n n^{\theta-1}$ for some constant $c_1>0$.
  By Theorem \ref{T3.1}(3) with $\dd=1$ and  $\ll(t)=\ss{\log\log t}$ for large $t>0$, we obtain
 $$
 \P\Big(S_{t_n}^\nu\ge (1+2r)\ss{2t_n\log\log t_n}\Big) \le c_2\exp\big[-(1+r)^2\log\log t_n\big]= c_2n^{-\theta(1+r)^2},\ \ n\ge 2
 $$
 for some constant $c_2>0$. Since $\theta (1+r)^2 >1$, this implies
 $$\sum_{n=2}^\infty \P\big(S_{t_n}^\nu\ge (1+2r)\ss{2t_n\log\log t_n}\big)<\infty,$$
 so that by the Borel-Cantelli lemma,
$$ \limsup_{n\to\infty} \ff{S_{t_n}^\nu}{\ss{2t_n\log\log t_n}}\le 1+2r,\ \ \text{a.s.}$$
Combining this with \eqref{3*} we obtain
$$ \limsup_{t\to\infty} \ff{S_{t}^\nu}{\ss{2t\log\log t}}\le 1+2r,\ \ \text{a.s.}$$
Since $r>0$ is arbitrary, we prove \eqref{U4}.

It remains to prove the following lower bound estimate for $\dd=1$ under condition \eqref{CV}:
 \beq\label{LB} \limsup_{t\to\infty} \ff{S_t^\nu}{\ss{2t\log\log t}}\ge  1,\ \ \text{a.s.}\end{equation}
 For any $\vv\in (0,\ff 1 2)$ and $p\ge 2$,   we have
 \beg{align*} &\lim_{l\to\infty} \ff{(1-2\vv)\ss{2p^l\log\log p^l}-(1-\vv)\ss{2(p^l-p^{l-1}-l)\log\log (p^l-p^{l-1}-l)}}{\ss{2(p^{l-1}+l)\log\log (p^{l-1}+l)}}\\
 &= (1-2\vv)\ss p -(1-\vv)\ss{p-1},\end{align*}
 which goes to $-\infty$ as $p\to\infty$. Then we may find constants $p,l_0\ge 2$ such that
 \beq\label{**0}\beg{split} &(1-2\vv)\ss{2p^l\log\log p^l}-(1-\vv)\ss{2(p^l-p^{l-1}-l)\log\log (p^l-p^{l-1}-l)}\\
  &\le -2 \ss{2(p^{l-1}+l)\log\log(p^{l-1}+l)},\ \ l\ge l_0.\end{split} \end{equation}
 Let
 $$G_l=\big\{S_{p^l}^\nu\ge (1-2\vv) \ss{2p^l\log\log p^l}\big\},\ \ l\ge l_0.$$ We aim to prove
 \beq\label{**1} \P\Big(\bigcap_{n=l_0}^\infty\bigcup_{m=n}^\infty G_m\Big)=1,\end{equation} which implies
 $$\limsup_{n\to\infty} \ff{S_{p^n}^\nu}{\ss{2p^n\log\log p^n}}\ge 1-2\vv,$$ so that
 $$\limsup_{t\to\infty} \ff{S_{t}^\nu}{\ss{2t\log\log t}}\ge 1-2\vv.$$ By the arbitrariness of $\vv\in (0,\ff 1 2)$ this implies the desired estimate \eqref{LB}.

 We now prove \eqref{**1}. For any $l\ge l_0+1,$ by \eqref{**0} we have
 $$G_l^c   \subset H_{l,1} \cup H_{l,2}$$
for
\beq\label{P1}\beg{split}
& H_{l,1}:= \Big\{S_{p^{l-1}+l}^\nu\le -2 \ss{2(p^{l-1}+l)\log\log (p^{l-1}+l)}\Big\}, \\
& H_{l,2}:=\Big\{S_{p^l}^\nu -S_{p^{l-1}+l}^\nu \le (1-\vv)\ss{2(p^l- p^{l-1}-l)\log\log (p^l-p^{l-1}-l)}\Big\}.
\end{split} \end{equation}
Hence, for any integer numbers $l>n\ge l_0$, by Markov property we have
 \beq\label{**2}
 \beg{split}
 \P\Big(\bigcap_{m=n}^l  G_m^c\Big) & \le \P \bigg(\left\{H_{l,1} \bigcup H_{l,2}\right\}\bigcap \left\{\bigcap_{m=n}^{l-1} G_m^c\right\} \bigg) \\
 & \le \P\Big(H_{l,1}\Big)+  \E\bigg[\P\Big(H_{l,2} \Big|\F_{p^{l-1}}\Big) \prod_{m=n}^{l-1}1_{G_m^c} \bigg]\\
 &=\P\Big(H_{l,1}\Big)+  \E\bigg[\P\Big(H_{l,2} \Big|X^\nu_{p^{l-1}}\Big) \prod_{m=n}^{l-1}1_{G_m^c} \bigg].
\end{split}\end{equation}
By Theorem \ref{T3.1}(3) with $\ll(t)= \ss{\log\log t}$ for large $t>0$, there exist  constants
 $c_1,c_2>0$ such that for $l \ge l_0$,
 \
 \beq\label{P2}\beg{split} &\P\big(H_{l,1}\big)\le c_1 \exp\big[-3\log\log (p^{l-1}+l)]\le c_2 l^{-3}.\end{split}\end{equation}
 Let $P_t(x,\cdot)$ be the distribution of $X_t^x$.
 \
 \beq\label{**2-1}
 \beg{split}
 & \ \ \ \ \E\bigg[\P\big(H_{l,2} \big|X^\nu_{p^{l-1}}\big) \prod_{m=n}^{l-1}1_{G_m^c}\bigg]\\
 &=\E\bigg[  \P\bigg(\ff{S_{p^l}^\nu -S_{p^{l-1}+l}^\nu} {\ss{2(p^l- p^{l-1}-l)\log\log (p^l-p^{l-1}-l)}}\le 1-\vv \bigg| X^\nu_{p^{l-1}}  \bigg)\prod_{m=n}^{l-1}1_{G_m^c}\bigg]\\
 &\le \P\big(|X_{p^{l-1}}^\nu|\ge l^{r_0}\big)\\
 &\quad + \E \bigg[\prod_{m=n}^{l-1}1_{G_m^c}\bigg]
   \sup_{|x|\le l^{r_0}}
 \P\bigg(\ff{S_{p^l-p^{l-1}-l}^{P_l(x,\cdot)} }{\ss{2(p^l- p^{l-1}-l)\log\log (p^l-p^{l-1}-l)}}  \le  1-\vv  \bigg),\end{split}\end{equation}
where    the last step follows from the  time-homogenous Markov property   that given $X_{p^{l-1}}^\nu=x$, the conditional distribution of $X_{p^{l-1}+l}^\nu$ is $P_l(x,\cdot)$, and the conditional distribution of $S_{p^l}^\nu -S_{p^{l-1}+l}^\nu$ coincides with the distribution of $S_{p^l-p^{l-1}-l}^{P_l(x,\cdot)}$.
By the condition in (4) we have
\
\beq\label{**3}
\beg{split}
\P\big(|X_{p^{l-1}}^\nu| \ge l^{r_0}\big)& \le l^{-2} \E |X_{p^{l-1}}^\nu|^{2/r_0} \\
&= \ff 1 {l^2} \int_{\R^d} \rr_\nu(x) P_{p^{l-1}}|\cdot|^{2/r_0}(x) \mu(\d x) \\
&\le \ff 1 {l^2} \left[\mu(\rr_\nu^q)\right]^{\ff 1 q} \left[\mu(|\cdot|^{\ff{2q}{r_0(q-1)}})\right]^{\ff{q-1}q} \le \frac{c_3}{l^2}
\end{split}
\end{equation} for some constant $c_3>0$.
Moreover,   by \eqref{CV} and Theorem \ref{T3.1}(3) with $\ll(t)=\log\log t$ for large $t$, we have
\beg{align*} &\sup_{|x|\le l^{r_0}} \P\Big(S_{p^l-p^{l-1}-l}^{P_l(x,\cdot)}   \le (1-\vv)\ss{2(p^l- p^{l-1}-l)\log\log (p^l-p^{l-1}-l)} \Big)\\
&= 1 - \inf_{|x|\le l^{r_0}}
 \P\Big(S_{p^l-p^{l-1}-l}^{P_l(x,\cdot)}   > (1-\vv)\ss{2(p^l- p^{l-1}-l)\log\log (p^l-p^{l-1}-l)} \Big)\\
 &\le 1-\exp\Big[-\Big(1-\ff \vv 2\Big)^2\log\log (p^l-p^{l-1}-l)\Big]\le 1-(l\log p)^{-(1-\ff\vv 2)^2},\ \ l\ge l_0+1. \end{align*}
 Combining this with \eqref{**2-1} and \eqref{**3}, we get
 \
  \beq\label{**2-2}
 \E\left[\P\Big(H_{l,2} \Big|X^\nu_{p^{l-1}}\Big) \prod_{m=n}^{l-1}1_{G_m^c}\right] \le \frac{c_3}{l^2}+\left(1-(l\log p)^{-(1-\ff\vv 2)^2} \right)\E \bigg(\prod_{m=n}^{l-1}1_{G_m^c}\bigg),
\end{equation}
which, together with \eqref{**2} and \eqref{P2}, yields
\
 \begin{equation*}
 \beg{split}
 \P\Big(\bigcap_{m=n}^l  G_m^c\Big) & \le \frac{c}{l^2}+\big(1-(l\log p)^{-(1-\ff\vv 2)^2}\big)\P\Big(\bigcap_{m=n}^{l-1}  G_m^c\Big), \ \ l\ge n+1\ge l_0+1
\end{split}\end{equation*}
 for some constants $c>0$. Therefore, by induction we obtain
 $$ \P\Big(\bigcap_{m=n}^l G_m^c\Big) \le c\sum_{m=n}^l m^{-2} +\prod_{m=n+1}^l  \big(1-(m\log p)^{-(1-\ff\vv 2)^2}),\ \ l\ge n+1\ge l_0+1.$$ This implies
 $$\P\Big(\bigcap_{m=n}^\infty G_m^c\Big) = \lim_{l\to\infty} \P\Big(\bigcap_{m=n}^l G_m^c\Big) \le c \sum_{m=n}^\infty m^{-2},$$ so that
 $$ \P\Big(\bigcap_{n=1}^\infty\bigcup_{m=n}^\infty G_m\Big)= \lim_{n\to\infty}  \P\Big( \bigcup_{m=n}^\infty G_m\Big) = 1-\lim_{n\to\infty} \P\Big(\bigcap_{m=n}^\infty G_m^c\Big)=1.$$
 Therefore, \eqref{**1} holds and  the proof is finished.
 \end{proof}

\beg{proof}[Proofs of Theorem  $\ref{T1.1}$ and Theorem $\ref{C1.2}$]  Let $\psi= 2\si^{-1}B-\si^*\nn\log\rr.$ By Proposition \ref{P2.1}, all conditions in Theorem \ref{T3.1} hold, so that Theorem \ref{T1.1} follows immediately. To prove Theorem \ref{C1.2}, we first show that  \beq\label{LT} \E \e^{r |S_{t_0}^x|}\le H(x),\ \ x\in\R^d \end{equation} holds for some constants $r,t_0>0$ and    locally bounded function $H$ on $\R^d$. Indeed, writing $\e^{r|S_{t_0}^x|}\le \e^{rS_{t_0}^x}+ \e^{-rS_{t_0}^x},$ \eqref{LT} with $r=1$ and small $t_0>0$ follows from   Proposition \ref{P2.1}(1) and (4) since
$$\E\e^{\int_0^{t_0}|\psi(X_s^x)|^2\d s} \le \ff 1 {t_0} \int_0^{t_0} \E \e^{t_0 |\psi(X_s^x)|^2}\d s =\ff 1 {t_0} \int_0^{t_0} P_s \e^{t_0 |\psi|^2}(x)\d s.$$
Now, by the Markov property, $(S_t^x-S_{t_0}^x)_{t\ge t_0}$    identifies with $(S_{t-t_0}^{P_{t_0}(x,\cdot)})_{t\ge t_0}$ in distribution. Moreover,
by Proposition \ref{P2.1},  $\mu(p_{t_0}(x,\cdot)^2)$  is locally bounded in $x$. So, (1) and (3) in Theorem \ref{C1.2}  follow from (2) and (4) in Theorem \ref{T3.1}  respectively. To prove Theorem \ref{C1.2}(2), we observe that for  $\ll$ therein it follows from \eqref{LT} that
$$\limsup_{t\to\infty} \ff 1{\ll(t)^2} \log \P\bigg(\ff{|S_{t_0}^x|}{\ll(t)\ss{2t}}>\vv\bigg) \le \limsup_{t\to\infty} \ff{\log \exp[-r \ll(t)\ss{2t}]}{\ll(t)^2}
 =-\infty$$ locally uniformly in $x$. Then Theorem \ref{C1.2}(2) follows from Theorem \ref{T3.1}(3) since $\mu(p_{t_0}(x,\cdot)^2)$ is locally bounded in $x$.

 \end{proof}

\section{SDEs with multiplicative noise}
 Consider the SDE
 \beq\label{E'} \d X_t = B(X_t)\d t+\si(X_t)\d W_t\end{equation}
 as \eqref{SDE}, but $\si$ now depends on the space variable $x$. When $B$ is Lipschitz continuous and $\si\in C_b^1(\R^d; \R^d\otimes \R^d)$ such that $c_1I\le \si\si^*\le c_2I$ for some constants $c_2\ge c_1>0$, and  the dissipativity condition
 \beq\label{C2} \|\si(x)-\si(y)\|^2_{HS} +2\<B(x)-B(y), x-y\> \le   -K|x-y|^2,\ \ x,y\in \R^d\end{equation} holds for some constant $K>0$,  we can prove the existence and uniqueness of invariant probability measure $\mu$ with $\mu(\e^{\vv|\cdot|^2})<\infty$ for some $\vv>0$. Moreover,   by \cite[Theorem 1.1]{Wan11}   the associated Markov semigroup $P_t$ satisfies
 \beq\label{HI'} |P_t f(x)|^p\le (P_t |f|^p(y))\exp\bigg[\ff{c(p) |x-y|^2}{\e^{Kt}-1}\bigg],\ \ p> p_0, t>0, x,y\in \R^d\end{equation} for all $f\in \B_b(\R^d)$, where $p_0\ge 1$ is a constant
 depending on  the smallest and largest eigenvalues of $\si\si^*$ and $c(p)$ is a constant depending on $p$.  As in the proof of Proposition \ref{P2.1}, this together with $\mu(\e^{\vv |\cdot|^2})<\infty$ implies the hyperboundedness   of $P_t$ as well as the local boundedness of $\mu(p_t(x,\cdot)^p)$ for some $p>1$.
 The only problem for us to extend Theorem  \ref{T1.1} and Theorem \ref{C1.2} to the present setting is   that we do not have good enough integration by parts formula to imply
 $\mu(\e^{\vv|\nn\log \rr|^2})<\infty$ for some constant $\vv>0$, where $\rr$ is the density of $\mu$ which is again strictly positive and belongs to $\cap_{p>1} W_{loc}^{p,1}(\d x)$ according to \cite{BR}.  Due to this problem, in the moment we are not able to start from a given drift $B$, but start from a given invariant probability measure $\mu$ with the required property
 $\mu(\e^{\vv|\nn\log \rr|^2})<\infty$. This can be done by perturbations to symmetric diffusion process.

 Now, let $V\in C^2(\R^d)$ such that
 \beq\label{V} \int_{\R^d} \e^{\vv |\nn V(x)|^2+V(x)}\d x<\infty\ \text{for\ some\ constant\ }\vv>0.\end{equation} Without loss of generality, we may and do assume that $\mu(\d x):= \e^{V(x)}\d x$ is a probability measure. Let $b: \R^d\to \R^d$ be locally integrable such that
 \beq\label{BB} \int_{\R^d}\<b,\nn f\>\d\mu=0,\ \ f\in C_0^\infty(\R^d).\end{equation}
 Let $\{e_i\}_{i=1}^d$ be the standard ONB of $\R^d$. We assume that
 \beq\label{BB2} B:= b+ \ff 1 2 \sum_{ij=1}^d \{\pp_j (\si\si^*)_{ij}\} e_i + \ff 1 2 (\si\si^*)\nn V\end{equation}
 has linear growth. Then $\mu$ is the unique invariant probability measure of $P_t$ and, as in the additive noise case,  the sample EPR of the solution to \eqref{E'} can be formulated as
$$\ep_t(X_{[0,t]}) = \ff 1 t \int_0^t \<\psi(X_s),\d W_s\> +\ff 1 {2t}\int_0^t |\psi(X_s)|^2\d s,$$ where, noting that $\nn\rr=V$ in the present setting,
$$\psi:=  2\si^{-1}B-\si\si^*\nn V=  2\si^{-1}b +\sum_{i,j=1}^d \{\pp_j (\si\si^*)_{ij}\} \si^{-1} e_i.$$
Since $B$ has linear growth, $c_1I\le \si\si^*\le c_2I$ and $\mu(\e^{\vv (|\cdot|^2+|\nn V|^2)})<\infty $ for some constant $\vv>0$, we have $\mu(\e^{\vv |\psi|^2})<\infty$ for some
$\vv>0$ as required in Theorem \ref{T3.1}. Therefore,    the following result follows from Theorem \ref{T3.1}.

 \beg{thm}\label{T4.1} Let $\si\in C_b^1(\R^d; \R^d\otimes \R^d)$ be such that $c_1I\le \si\si^*\le c_2I$ for some constants $c_2\ge c_1>0$, let $V\in C^2(\R^d)$ satisfy
 $\eqref{V}$,   and let $b$ satisfy $\eqref{BB}$ such that $B$ defined in $\eqref{BB2}$ has linear growth.
If $\eqref{C2}$ holds for some constant $K>0,$ then all assertions in Theorems $\ref{T1.1}$ and $\ref{C1.2}$ hold for $\log\rr=V$. \end{thm}

\paragraph{Example 4.1.} A simple example satisfying all conditions in Theorem \ref{T4.1} is that $V(x)= \aa-\bb |x|^2$ for constants $\aa\in\R$ and $\bb>0$ such that $\e^{V(x)}\d x$ is a probability measure, $b(x) =A x$ for some antisymmetric $d\times d$-matrix $A$,   $\si\in C_b^1(\R^d; \R^d\otimes \R^d)$ such that  $c_1I\le \si\si^*\le c_2I$ for some constants $c_2\ge c_1>0$, and
\beq\label{LTZ} \ll < \ff 1 {2c_2}\bigg( d \|\nn\si\|_{\infty}^2 + 2\|A\| + \Big\{\sum_{i=1}\Big(\sum_{j=1}^d \|\pp_j (\si\si^*)_{ij}\|_\infty \Big)^2\Big\}^{\ff 1 2}\bigg), \end{equation}
  which implies \eqref{C2} for some constant $K>0$.

  \

Below we use the log-Sobolev inequality to replace the condition \eqref{C2}.
As observed in the proof of Proposition \ref{P2.1}(3)  that if $P_t$ has an invariant probability measure $\mu$ then $1$ is a simple eigenvalue of $P_t$ for $t>0$. So, the key condition in Theorem \ref{T3.1} is the hyperboundedness of $P_t$. We will see that this follows from the following log-Sobolev inequality in the uniformly elliptic  case:
\beq\label{LS0} \mu(f^2\log f^2)\le C \mu(|\nn f|^2),\  \ f\in C_0^\infty(\R^d), \mu(f^2)=1\end{equation} for a constant $C>0.$
Due to the Bakry-Emery citerion \cite{BE}, this inequality holds with $C=\ff 2 K$ if $\Hess_V\le -K$ for some constant $K>0$. By \cite{CW} the log-Sobolev inequality
\eqref{LS0} holds if $\Hess_V(x) \le -K$ for some constant $K>0$ and large enough $|x|>0$. In case that $\Hess_V(x)\le K$ for some positive constant $K>0$ and large enough $|x|$, according to \cite{W01} the log-Sobolev inequality holds provided $\mu(\e^{\vv |\cdot|^2})<\infty $ holds for some $\ll>\ff K 2.$ See also \cite{CGW} for Lyapunov type sufficient conditions of the log-Sobolev inequality.

Now, we state the following alternative version of Theorem \ref{T4.1}   with condition \eqref{C2} replaced by \eqref{LS0}. This result   applies to Example 4.1 without assuming   \eqref{LTZ}.  The price we have to pay is that we  can not prove the exact LIL due to the lack of  the moment estimate \eqref{HEAT} on the heat kernel.

\beg{thm}\label{T4.2} Let $\si\in C^1_b(\R^d; \R^d\otimes \R^d)$ such that  $c_1I\le \si\si^*\le c_2I$ for some constants $c_2\ge c_1>0$,   let $V\in C^1(\R^d)$ such that $\mu(\d x):=\e^{V(x)}\d x$ is a probability measure satisfying $\eqref{V}$ and $\eqref{LS0}$,    and let $b$  satisfy $\eqref{BB}$ and has linear growth.  Then for the $SDE$ with $B$ given by \eqref{BB2}, all assertions in Theorems $\ref{T1.1}$ (1)-(3) hold, and  for any $\nu\in\scr P(\R^d)$ with $\rr_\nu:=\ff{\d\nu}{\d\mu}\in L^q(\mu)$ for some $q>1$,   $\P$-a.s.
$$  \limsup_{t\to\infty} \ff{S_t^\nu}{\ss{2t\log\log t}}  \le\ss{\dd},\ \   \liminf_{t\to\infty} \ff{S_t^\nu}{\ss{2t\log\log t}} \ge-\ss\dd.$$  \end{thm}

\beg{proof} Since $\si$ is $C^1$-smooth and $B$ has linear growth, the SDE \eqref{E'} has a unique non-explosive solution. Let $P_t$ be the associated Markov semigroup. By \eqref{LS0},    $\mu(\e^{\vv|\cdot|^2})<\infty$ holds for some constant $\vv >0$ (see \cite{AMS}).  Since  $\psi:= 2\si^{-1}B- \si^* \nn V$, $B$ has linear growth and $c_1 I\le \si\si^*\le c_2 I$, this and \eqref{V} imply   $\mu(\e^{\vv|\psi|^2})<\infty$   for some $\vv>0$.

Moreover,  as explained in the proof of Proposition \ref{P2.1}(3) using \cite[Proposition 2.11]{MR},   \eqref{BB}, \eqref{BB2} and $b\in L^p(\mu)$ for all $p>1$ implies that $P_t$ is associated to a Dirichlet form with symmetric part
$$ \EE(f,g):= \mu(\<\si\si^*\nn f,\nn g\>), \ \ f,g\in H_\si^{2,1}(\mu),$$ and $1$ is a simple eigenvalue of $P_t$ for $t>0.$
By \eqref{LS0}   we have
\beq\label{LSI} \mu(f^2\log f^2) \le \ff {C}{c_1}\EE(f,f),\ \ f\in H_\si^{2,1}(\mu), \mu(f^2)=1.\end{equation}
 So, according to \cite{Gross76}, the semigroup $P_t$ is hypercontractive.
 Then the proof is sinished by Theorem \ref{T3.1}. \end{proof}

\end{document}